\documentclass[english,a4paper,10pt]{article}

\usepackage{xcolor}
\usepackage{amsfonts,amsmath,amstext,amsbsy,amsthm,amscd}
\usepackage{graphics}
\usepackage{epsfig}
\usepackage{geometry}
\usepackage{amssymb}
\usepackage[latin1]{inputenc}
\usepackage{babel,indentfirst,amsfonts,array}

\newtheorem{thm}{Theorem}[section]
\newtheorem{lemme}[thm]{Lemma}
\newtheorem{prop}[thm]{Proposition}
\newtheorem{cor}[thm]{Corollary}

\def\P{\mathbb P}
\def\E{\mathbb E}

\def\R{\mathbb R}

\def\Z{\mathbb Z}
\def\N{\mathbb N}
\def\({\left(}
\def\){\right)}
\def\[{\left[}
\def\]{\right]}
\def\fin{\hfill\square}

\def\fin{\hfill $\square$}

\title{\textsc{On homogeneous and oscillating random walks on the integers}
\author{Julien Br\'emont}
\date{Universit\'e Paris-Est Cr\'eteil,~janvier 2022, version 2}
}

\begin{document}

\maketitle

\setcounter{page}{1}

\begin{abstract}
We study the recurrence of homogeneous and oscillating random walks on the integers, simplifying former works of Spitzer and Kemperman, respectively. We add general remarks and discuss some links with renewal theory.

\end{abstract}

\footnote{
\begin{tabular}{l}\textit{AMS $2020$ subject classifications~: 60G17, 60J10.} \\
\textit{Key words and phrases~: Homogeneous random walk, oscillating random walk, recurrence criterion.} 
\end{tabular}}

\section{Introduction}

A primary question in the study of the asymptotic behaviour of a discrete time Markov chain on $\Z^d$ is that of recurrence. When the jumping law is the same everywhere (homogeneous case), this problem concerns Birkhoff sums $S_n=X_1+\cdots+X_n$, where the $(X_i)$ are $\Z^d$-valued random variables, independent and identically distributed ($i.i.d.$), with common law $\mu$. Improving a former result of Chung and Fuchs, Spitzer's analytical recurrence criterion (1957, cf \cite{spitz}, T2) states that transience is equivalent to the integrability of $Re(1/(1-\hat{\mu}))$ on the unit cube of $\R^d$. Importantly, this result doesn't require any moment condition. For a general model of inhomogeneous random walk on $\Z^d$, the transition law at $x\in\Z^d$ is given by some probability measure $\mu_x$. The question of the recurrence is very delicate, as one has to understand how the Markov chain perceives the space (or ``environment") given by the couple ($\Z^d$, $(\mu_x)_{x\in\Z^d}$).

\medskip
A situation where the $i.i.d.$ case helps is the following one. Fix a general Markov chain on $\Z^d$, starting at $0$, and begin with a discussion of a naive recursive approach of the problem. Given subsets $\Z^d=F_0\supset\cdots\supset F_K=\{0\}$, the random walk is recurrent at $0$ if and only if, inductively on $0\leq k<K$, the induced random walk in $F_{k}$ visits $F_{k+1}$ infinitely often, almost-surely. A difficulty is that, quasi inevitably, the induced Markov chains are heavy-tailed. This strategy works for instance for random walks in a stratified environment, as considered in \cite{CP2003, jb1,jb2}, as there is in this case a natural filtration. A typical example is a nearest-neighbour Markov chain in $\Z^2$, when the transition laws only depend on the second coordinate. In this situation the vertical component of the random walk, in restriction to vertical movements, is a nearest-neighbour one-dimensional Markov chain. The necessary and sufficient condition for its recurrence is well-known, for example in the theory of birth and death processes, and therefore corresponds to the recurrence of the initial random walk in $\Z\times\{0\}$ (we are thus led to taking $\Z^2=F_0\supset F_1=\Z\times\{0\}\supset F_2=\{0\}$). When the condition holds, the induced random walk in $\Z\times\{0\}$ is heavy-tailed, but $i.i.d.$, due to the invariance of the environment by horizontal translations. The analysis is then naturally orientated into precising the jumping law of this random walk and applying Spitzer's recurrence criterion for $\Z$-valued $i.i.d.$ sums; cf \cite{jb1,jb2}. 

\medskip
A point that is not a detail is that the proofs of \cite{jb1,jb2} would have been seriously more delicate to handle if we had to use the Chung-Fuchs result in place of Spitzer's criterion, and this orientates the present work. To now make a small step outside stratified random walks in the plane, one has to develop results around the cornerstone that constitutes Spitzer's theorem. As a first natural extension, it appears important to prove for Kemperman's oscillating random walk \cite{kemperman} a result in the same spirit. Kemperman's results \cite{kemperman} on this model indeed correspond to the Chung-Fuchs theorem for homogeneous random walks. 

\medskip
A preliminary step of clarification is necessary concerning Spitzer's theorem. The known proof, available in \cite{spitz}, is not linear and requires to enter to some extent the apparatus of Potential Theory for discrete Markov chains. It is also disseminated in \cite{spitz} and a little painful to reconstitute (Kesten-Spitzer \cite{kest-spitz} helps). We provide here a short(er) proof of Spitzer's result and highlight that it consists in computing some kind of second derivative at infinity of the Green function in two different ways, a probabilistic one and one relevant from Fourier Analysis. Invigorating a lemma due to Chung from the Potential Theory of discrete recurrent Markov chains, we show that the probabilistic part of the proof is in fact very general. The harmonic part will not be discussed here. We next point out some links with renewal theory and, in the last section, consider oscillating random walks on the integers, reproving the results of Kemperman \cite{kemperman} on the recurrence of this model. A combinatorial remark allows a shortcut in the analysis. 

\medskip
This text is a mainly revisit, the material and the results being essentially not new. Our effort has been concentrated on the exposition, which tries to be in straight line and self-contained. Many questions are addressed along the way, essentially on inner products of probabilistic Green functions and their translations in Harmonic analysis.

\medskip
We fix $\Z$ as state space, except for section \ref{probgen}. We now recall classical facts and notations.

\section{Preliminaries}
{\it 1) Laws and characteristic functions.} Let $X$ be a $\Z$-valued random variable, with law written as ${\cal L}(X)=\mu$, and $(X_n)_{n\geq1}$ be $i.i.d.$ copies. We suppose $X$ non constant, with $\mbox{gcd}(\mbox{Supp}(\mu))=1$. Let the characteristic function $\hat{\mu}(t)=\E(e^{itX}),~t\in\R$, a $2\pi$-periodic function. Defining the integer $d=\mbox{gcd}(\mbox{Supp}({\cal L}(X_1-X_2)))\geq1$, then~:

$$\hat{\mu}(t)=1\mbox{ iff }t\in 2\pi\Z\mbox{ and }|\hat{\mu}(t)|=1\mbox{ iff }t\in(2\pi/d)\Z.$$

\smallskip
\noindent
Write $S_n=\sum_{i=1}^nX_i$, with $S_0=0$. The property that $|\hat{\mu}(t)|<1$ on $[0,2\pi]$, except for finitely many $t$, implies that ${\cal L}(S_n)$ does not concentrate around any point, as $n\rightarrow+\infty$~:

$$\forall y\in\Z,~P(S_n=y)=\frac{1}{2\pi}\int_0^{2\pi}e^{-ity}(\hat{\mu}(t))^n~dt\rightarrow_{n\rightarrow+\infty}0.$$

\noindent
Also, for some $\alpha>0$ and small $t$, $Re(1-\hat{\mu}(t))\geq \alpha t^2$. Indeed, take $M>0$ so that $P(0<|X|<M)>0$. Then for $|t|\leq\pi/M$, $Re(1-\hat{\mu}(t))= 2\E(\sin^2(tX/2))\geq 2\pi^{-2}t^2\E(X^21_{|X|<M})$.

\bigskip
{\it 2) Markov chains and Green functions.} For any Markov chain $(S_n)$ on $\Z$, $P_x$ and $\E_x$ stand for $x$ as starting point. The Green function is $G(x,y)=\E_x(\sum_{n\geq0}1_{S_n=y})$. Let also $G_N(x,y)=\E_x(\sum_{0\leq n<N}1_{S_n=y})$, $N\geq1$. For $y\in\Z$, set $T_y=\min\{n\geq1,~S_n=y\}$. For any $x\not=y$~:

\begin{equation}
\label{deb}
G_N(x,y)=\sum_{1\leq k< N}P_x(T_y=k)G_{N-k}(y,y).
\end{equation}

\noindent
Hence, $G_N(x,y)\leq G_N(y,y)$ and $G(x,y)=P_x(T_y<\infty)G(y,y)$. Also, for $x\in\Z$~:

$$G(x,x)=\sum_{n\geq0}P_x(T_x<\infty)^n=1/(1-P_x(T_x<\infty)),$$

\noindent
so the recurrence of $x$, i.e. the property $P_x(T_x<\infty)=1$, is equivalent to $G(x,x)=+\infty$. For any $x\not=y$ with $P_x(T_y<\infty)>0$, note that $P_x(T_x<T_y)<1$. Then, in the same way~:

\begin{equation}
\label{troncat}
\E_x\({\sum_{n=0}^{T_y-1}1_{S_n=x}}\)=1+\sum_{n\geq1}P_x(T_x<T_y)^n=\frac{1}{1-P_x(T_x<T_y)}<\infty.
\end{equation}

\noindent
Still for any $x\not=y$, we have~:

\begin{eqnarray}
\label{import}
G_N(x,x)&=&\E_x\({\sum_{n=0}^{T_y\wedge N-1}1_{S_n=x}}\)+\E_x\({1_{T_y<N}\sum_{n=T_y}^{N-1}1_{S_n=x}}\)\nonumber\\
&=&\E_x\({\sum_{n=0}^{T_y\wedge N-1}1_{S_n=x}}\)+\sum_{k=1}^{N-1}P_x(T_y=k)G_{N-k}(y,x).
\end{eqnarray}

\noindent
Thus, $0\leq G_N(x,x)-G_N(y,x)\leq \E_x(\sum_{n=x}^{T_y\wedge N-1}1_{S_n=x})$. We deduce the important {\it claim~:} for $x\not=y$ with $P_x(T_y<\infty)>0$, then $(G_N(x,x)-G_N(y,x))_{N\geq0}$ is bounded.

\medskip
In the particular case when the chain is homogeneous, $G(x,y)=G(x-y,0)$ and $G_N(x,y)=G_N(x-y,0)$. Notice also that for $x\not=0$, we have $P_0(T_0<T_x)=P_0(T_0<T_{-x})$, as~:

\begin{eqnarray}
P_0(T_0<T_x)&=&\sum_{k\geq1}P_0(S_k=0, S_l\not\in\{0,x\}, 0<l<k)\nonumber\\
&=&\sum_{k\geq1}P_0(S_k=0, S_k-S_l\not\in\{0,-x\}, 0<l<k)=P_0(T_0<T_{-x}).\nonumber
\end{eqnarray}

\noindent
Still in the homogeneous case, with a step $X$ of law $\mu$, we often put $\mu$ as a superscript and write $S_n^{\mu}$, $G^{\mu}(x,y)$, $G_N^{\mu}(x,y)$, as well as $\E^{\mu}(f(X))$ for $\int_{\Z}fd\mu$.

\section{Homogeneous case : Spitzer's analytical criterion}

Let $(S_n)$ be a homogeneous random walk on $\Z$ with step $\mu$, not Dirac and $\mbox{gcd}(\mbox{Supp}(\mu))=1$. On $(0,2\pi)$, the function $t\longmapsto\mbox{Re}\({1/(1-\hat{\mu}(t))}\)$ is $>0$, continuous and symmetric under $t\longmapsto 2\pi-t$. It belongs to $L^1(0,2\pi)$ iff it is in $L^1(0,\varepsilon)$, for some $\varepsilon>0$.

\begin{thm}(Spitzer, 1957)

\noindent
The point $0$ is transient for $(S_n)$ iff $\int_0^{2\pi}\mbox{Re}\({1/(1-\hat{\mu}(t))}\)dt<+\infty$.\end{thm}

\medskip
This will follow from the next proposition, where constants are optimal.

\begin{prop} 
\label{inega}
We have $\displaystyle G(0,0)\leq \frac{1}{\pi}\int_0^{2\pi}\mbox{Re}\({1/(1-\hat{\mu}(t))}\)dt\leq 2G(0,0)$.\end{prop}

\medskip
\noindent
{\it Proof of the proposition~:} 

\noindent
Set $a_N(y)=G_N(0,0)-G_N(0,y)=G_N(0,0)-G_N(-y,0)$. We show that $(a_N(y))_{N\geq0}$ is bounded. Take $y\not=0$. This is clear if $0$ is transient and if it is recurrent, then $P_0(T_{-y}<\infty)=1$, as $\mbox{gcd}(\mbox{Supp}(\mu))=1$. Thus $(a_N(y))_{N\geq0}$ is bounded, by the {\it claim} above.

\medskip
\noindent
{\it Step 1.} Let $x>0$. We show that $\Delta(x):=\lim_{N\rightarrow+\infty}(a_N(x)+a_N(-x))$ exists. We have~:

$$a_N(x)+a_N(-x)=\frac{1}{2\pi}\int_0^{2\pi}(2-e^{-itx}-e^{itx})\sum_{n=0}^{N-1}(\hat{\mu}(t))^ndt=\frac{1}{\pi}\int_0^{2\pi}\frac{1-\cos (tx)}{1-\hat{\mu}(t)}(1-(\hat{\mu}(t))^{N})dt.$$

\noindent
Since $|1-\hat{\mu}(t)|\geq Re(1-\hat{\mu}(t))\geq\alpha t^2$ and $x$ is fixed, $(1-\cos (tx))/(1-\hat{\mu}(t))$ is integrable. As $|\hat{\mu}(t)|<1$ except for finitely many values of $t$, the required limit exists and satisfies~:

\begin{equation}
\label{lim1}
\Delta(x)=\frac{1}{\pi}\int_0^{2\pi}\frac{1-\cos (tx)}{1-\hat{\mu}(t)}~dt=\frac{1}{\pi}\int_0^{2\pi}(1-\cos(tx))Re((1-\hat{\mu}(t))^{-1})dt.
\end{equation}

\noindent
{\it Step 2.} For $x>0$, we give a probabilistic expression for $\Delta(x)$. First, if $y\not=0$, using \eqref{deb} and \eqref{import}~:

$$a_N(y)=\E_0\({\sum_{n=0}^{T_y\wedge N-1}1_{S_n=0}}\)+\sum_{k=1}^{N-1}P_0(T_y=k)(G_{N-k}(y,0)-G_{N-k}(y,y)).$$

\noindent
By homogeneity, $a_N(y)=\E_0(\sum_{n=0}^{T_y\wedge N-1}1_{S_n=0})-\sum_{1\leq k<N}P_0(T_y=k)a_{N-k}(-y)$. Taking $y=x$ and adding $a_N(-x)$ we obtain~:

$$a_N(x)+a_N(-x)=\E_0\({\sum_{n=0}^{T_x\wedge N-1}1_{S_n=0}}\)+\sum_{k=1}^{N-1}P_0(T_x=k)(a_{N}(-x)-a_{N-k}(-x))+P_0(T_x\geq N)a_N(-x).$$

\noindent
Consider the terms on the right, when $N\rightarrow+\infty$. The first one tends to $\E_0(\sum_{n=0}^{T_x-1}1_{S_n=0})$. As $P_0(S_n=y)\rightarrow_{n\rightarrow+\infty}0$, we get $\lim_{N\rightarrow+\infty}a_{N+1}(y)-a_N(y)=0$ and thus $\lim_{N\rightarrow+\infty}a_{N}(-x)-a_{N-k}(-x)=0$ for fixed $k$. By dominated convergence, the second term goes to zero, as $(a_N(-x))_{N\geq 0} $ is bounded. The latter also implies that the third term goes to zero in case of recurrence and to $P_0(T_x=\infty)G(0,0)(1-P_0(T_{-x}<\infty))$ in case of transience. Thus, if $x>0$~:

\begin{equation}
\label{lim2}
\Delta(x)=\E_0\({\sum_{n=0}^{T_x-1}1_{S_n=0}}\)+1_{TR}G(0,0)P_0(T_x=\infty)P_0(T_{-x}=\infty).
\end{equation}

\noindent
{\it Step 3.} By \eqref{lim1}=\eqref{lim2}, for any $\delta>0$, $\pi^{-1}\int_\delta^{2\pi-\delta}(1-\cos(tx))Re((1-\hat{\mu}(t))^{-1})dt\leq2G(0,0)$. When $x\rightarrow+\infty$, we get $\pi^{-1}\int_{[\delta,2\pi-\delta]}Re((1-\hat{\mu}(t))^{-1})dt\leq 2G(0,0)$, by the Riemann-Lebesgue lemma. Letting $\delta\rightarrow0$, we get the second inequality. For the other direction, by \eqref{lim2}=\eqref{lim1}~:

$$\E_0\({\sum_{n=0}^{T_x-1}1_{S_n=0}}\)\leq(1/\pi)\int_0^{2\pi}(1-\cos tx)Re((1-\hat{\mu}(t))^{-1})dt.$$

\noindent
If $Re((1-\hat{\mu}(t))^{-1}\in L^1(0,2\pi)$ (if not, this is obvious), then letting $x\rightarrow+\infty$ and using again the Riemann-Lebesgue lemma for the right-hand side, we obtain the first inequality.  \fin

\bigskip
\noindent
\begin{remark} 
When transience holds, constants in Prop. \ref{inega} are optimal. If $\mbox{Supp}(\mu)\subset \N^*$, then $G(0,0)=1$ and $\Delta(x)=2-P_0(T_x<\infty)\rightarrow2-1/\E(X)$, as $x\rightarrow+\infty$, by renewal theory (this is reproved later in the paper). As $\lim_{x\rightarrow+\infty}\Delta(x)=\pi^{-1}\int_0^{2\pi}Re((1-\hat{\mu}(t))^{-1})dt$, the conclusion comes from the fact that $\E(X)$ can take any value in $[1,+\infty]$. 
\end{remark}

\medskip
\noindent
\begin{remark}
The idea, used by Spitzer, of approaching $G(0,0)$ in two steps, first by the finite $\lim_{N\rightarrow+\infty}(2G_N(0,0)-G_N(0,x)-G_N(0,-x))$ and next the limit as $x\rightarrow+\infty$, is classical and profound. A similar one is developed by Riemann in the first chapters of the theory of trigonometric series. This kind of ``second derivative at infinity'' suggests some general link between transience and a condition of positive curvature at infinity (where the terms are to be redefined). 
\end{remark}

\medskip
\noindent
\begin{remark}
The weak form of the theorem, due to Chung and Fuchs (1951), can be reduced to the following observation, where interversion is direct for $0<s<1$~: 

\begin{eqnarray}
\label{cf}
G^{\mu}(0,0)=\lim_{s\uparrow1}\sum_{n\geq0}s^nP^{\mu}(S_n=0)&=&\lim_{s\uparrow1}\sum_{n\geq0}\frac{1}{2\pi}\int_0^{2\pi}s^n(\hat{\mu}(t))^n~dt\nonumber\\
&=&\lim_{s\uparrow 1}\frac{1}{2\pi}\int_0^{2\pi}Re\({\frac{1}{1-s\hat{\mu}(t)}}\)dt.
\end{eqnarray}

\noindent
The finiteness of the right-hand side is thus a transience criterion for the random walk. Notice that the operation $s\uparrow1$ is not natural in this problem, as the level sets of $z\longmapsto Re(1/(1-z))$ in the unit disk are horocycles (Euclidean circles). There is no monotony in the limit and indeed, as seen above, the right-hand side may differ from $(2\pi)^{-1}\int_0^{2\pi}Re((1-\hat{\mu}(t))^{-1})dt$.
\end{remark}

\medskip
\noindent
\begin{remark} The theorem has been extended to general countable discrete Abelian groups by Kesten and Spitzer \cite{kest-spitz}, to $\R^d$ by Ornstein \cite{ornstein} and Port and Stone \cite{stone}. Recall also that in any case $\lim_{N\rightarrow+\infty}a_N(x)$ exists and is called the potential kernel; see Spitzer \cite{spitz}, chap. 7. \end{remark}
 
\medskip
\noindent
\begin{remark}In the second step of the proof of the proposition and in the transient case, one can directly write $\Delta(x)=G(0,0)(2-P_0(T_x<\infty)-P_0(T_{-x}<\infty))$, when $x>0$. It is interesting to check equality with \eqref{lim2} in this case, i.e. that for $x\not=0$~:

$$G(0,0)(2-P_0(T_x<\infty)-P_0(T_{-x}<\infty))=\frac{1}{1-P_0(T_0<T_x)}+G(0,0)P_0(T_x=\infty)P_0(T_{-x}=\infty).$$
 
 \noindent
Equivalently, for fixed $x\not=0$~:
 
 \begin{equation}
 \label{green}
 G(0,0)=\frac{1}{1-P_0(T_0<T_x)}\times \frac{1}{1-P_0(T_x<\infty)P_0(T_{-x}<\infty)}.
 \end{equation}
  
\noindent
This relation, valid for a general Markov chain (transient or not), corresponds to the following decomposition of $G(0,0)$. Let $T^{(0)}=0$ and next $T^{(k+1)}$ be the first time $>T^{(k)}$ of passage at zero after having visited $x$ at least once, for $k\geq0$. Then (the $k^{th}$ term being 0 if $T^{(k)}=+\infty$)~:

$$G(0,0)=\sum_{k\geq0}\E_0\({\sum_{n=T^{(k)}}^{T^{(k+1)}-1}1_{S_n=0}}\).$$

\medskip
\noindent
Call $A_k$ the generic term in the above sum. Then~:

$$A_0=\E_0\({\sum_{n=0}^{T^{(1)}-1}1_{S_n=0}}\)=\E_0\({\sum_{n=0}^{T_x-1}1_{S_n=0}}\)=\frac{1}{1-P_0(T_0<T_x)}.$$

\medskip
\noindent
Next, for $k\geq1$, $A_k=P_0(T^{(k)}<\infty)A_0=(P_0(T^{(1)}<\infty))^kA_0$. Now, $P_0(T^{(1)}<\infty)=P_0(T_x<\infty)P_0(T_{-x}<\infty)$ and this gives the announced formula when summing on $k\geq0$. Remark also that the last product is related with loops. More precisely~:
 
 $$P_0(T_{x}<\infty)P_0(T_{-x}<\infty)=\frac{P_0(T_x<T_0<\infty)}{1-P_0(T_0<T_x)}.$$
 
 \noindent
The left-hand side is $P_0(\mbox{reach }x\mbox{ and come back at }0)$. Decomposing it as the probability of making first $n\geq 0$ loops at 0 without touching $x$, then going directly to $x$ and finally coming back to 0, this equals $\sum_{n\geq0}P_0(T_0<T_x)^nP_0(T_x<T_0<\infty)$, so the right-hand side expression.
 
\medskip
Notice furthermore that one may readily derive from \eqref{green}, convening that $1/0^+=+\infty$, the following always valid form for $\Delta(x)$, $x>0$~:

$$\Delta(x)=\E_0\({\sum_{n=0}^{T_x-1}1_{S_n=0}}\)\({1+\({1+\frac{P_0(T_x<\infty)}{P_0(T_x=\infty)}+\frac{P_0(T_{-x}<\infty)}{P_0(T_{-x}=\infty)}}\)^{-1}}\).$$
 
\end{remark} 
 
\section{A general probabilistic result}
\label{probgen}
The proof of Proposition \ref{inega} consists in computing some second derivative of the Green function in two different ways, an analytical one, giving \eqref{lim1}, and a probabilistic one, leading to \eqref{lim2}. We show here that the probabilistic part is very general.

\medskip
Consider a general irreducible Markov chain $(S_k)_{k\geq0}$ on a countable state space. Fix two points $x\not=y$ and set~:

$$c_N=\E_x\({\sum_{n=0}^{T_y\wedge N-1}1_{S_n=x}}\)\mbox{ and }d_N=\E_y\({\sum_{n=0}^{T_x\wedge N-1}1_{S_n=y}}\).$$

\noindent
Then $c_N\uparrow c:=\E_x(\sum_{n=0}^{T_y-1}1_{S_n=x})$ and $d_N\uparrow d:=\E_y(\sum_{n=0}^{T_x-1}1_{S_n=y})$, finite quantities. Taking $N\geq3$, let us develop relation \eqref{import}, namely~:

\begin{eqnarray}
\label{larec}
G_N(x,x)&=&c_N+\sum_{k=1}^{N-1}P_x(T_y=k)G_{N-k}(y,x)\nonumber\\
&=&c_N+\sum_{k=1}^{N-1}P_x(T_y=k)\sum_{l=1}^{N-k-1}P_y(T_x=l)G_{N-k-l}(x,x)\nonumber\\
&=&c_N+\sum_{m=2}^{N-1}G_{N-m}(x,x)R_m,\end{eqnarray}

\noindent
where $R_m=\sum_{k,l\geq1, k+l=m}P_x(T_y=k)P_y(T_x=l)$, symmetric in $x$ and $y$. Notice that $\sum_{m\geq2}R_m=P_x(T_y<\infty)P_y(T_x<\infty)\leq1$, with equality iff the random walk is recurrent. In the same way~:

$$G_N(y,y)=d_N+\sum_{m=2}^{N-1}G_{N-m}(y,y)R_m.$$

\noindent
Let us show that $(dG_N(x,x)-cG_N(y,y))_{N\geq0 }$ is bounded. Set $u_N=dG_N(x,x)-cG_N(y,y)$ and $\varepsilon_N=dc_N-cd_N$. We shall prove that for some $C>0$, $\forall n\geq2$, $|\varepsilon_n|\leq C\sum_{m\geq n}R_m$.

\medskip
\noindent
Supposing this true, fix $N\geq3$ and maybe increase $C$ so that $|u_n|\leq C$, for $n<N$. The equality $u_N=\varepsilon_N+\sum_{2\leq m<N}R_mu_{N-m}$ then furnishes ~:

$$|u_N|\leq |\varepsilon_N|+\sum_{2\leq m<N}R_m|u_{N-m}|\leq C\sum_{m\geq N}R_m+\sum_{2\leq m<N}R_mC=C.$$

\noindent
The property $|u_n|\leq C$ is thus transmitted by recursion on $n\geq N$, giving the required boundedness. To establish the missing point, write $\varepsilon_N=d(c_N-c)-c(d_N-d)$ and note that~:

\begin{eqnarray}
c-c_N=\E_x\({1_{T_y>N}\sum_{k=N}^{T_y-1}1_{S_k=x}}\)
&=&\E_x\({1_{S_1\not=y,\cdots,S_N\not=y}\E_{S_N}\({\sum_{k=0}^{T_y-1}1_{S_k=x}}\)}\)\nonumber\\
&=&\E_x\({1_{T_y>N}P_{S_N}(T^*_x<T_y)c}\),\nonumber\end{eqnarray}

\noindent
with $T^*_x=\min\{n\geq0~|~S_n=x\}$. Hence $0\leq c-c_N\leq cP_x(T_y>N)$. We conclude with the remark that $\sum_{m\geq N}R_m\geq P_x(T_y>N)P_y(T_x>0)$, where $P_y(T_x>0)>0$.

\medskip
We obtain the following Doeblin type ratio limit theorem (cf Revuz \cite{revuz}, chap.4, ex. 4.10).

\begin{lemme}
\label{41}

$ $

\noindent
For any irreducible Markov chain on a countable state space and any points $x$ and $y$~:

$$\lim_{N\rightarrow+\infty}\frac{G_N(x,x)}{G_N(y,y)}=\alpha(x,y)\mbox{, with }\alpha(x,y):=\frac{\E_x(\sum_{0\leq n<T_y}1_{S_n=x})}{\E_y(\sum_{0\leq n<T_x}1_{S_n=y})}.$$

\smallskip
\noindent
Moreover, $\alpha(x,y)=G(x,x)/G(y,y)$ in the transient case and $\alpha(x,y)=\pi(x)/\pi(y)$ in the recurrent case, where $\pi$ is the unique (up to a positive multiple) invariant $\sigma$-finite measure.

\end{lemme}

\medskip
\noindent
{\it Proof of the lemma~:}

\noindent
In the transient case, directly from relation \eqref{larec}, $G_N(x,x)\rightarrow G(x,x)=c/(1-\sum_{m\geq2}R_m)$. Idem, $G_N(y,y)\rightarrow G(y,y)=d/(1-\sum_{m\geq2}R_m)$, giving the result. In the recurrent case, this follows from the boundedness of $(dG_N(x,x)-cG_N(y,y))$. In this situation, as $\pi$ is proportional to $z\longmapsto \E_y(\sum_{0\leq n<T_y}1_{S_n=z})$, we obtain~:

$$
\frac{\pi(x)}{\pi(y)}=\E_y\({\sum_{n=0}^{T_y-1}1_{S_n=x}}\)/1=P_y(T_x<T_y)\E_x\({\sum_{n=0}^{T_y-1}1_{S_n=x}}\)=\frac{1-P_y(T_y<T_x)}{1-P_x(T_x<T_y)}.$$

\noindent
Via \eqref{troncat}, we recognize $\alpha(x,y)$ and this concludes the proof of the lemma.\fin

\medskip
\noindent
\begin{remark}
When the regime is known (recurrence or transience), then $\alpha(x,y)$ has the form of a function of $x$ divided by the same function of $y$. This seems unclear a priori.
\end{remark}

\begin{lemme}
\label{fonda}

$ $

\noindent
Consider an irreducible and aperiodic Markov chain $(S_n)$ on a countable state space. Define $a_N(x,y)=G_N(x,x)-G_N(y,x)\geq0$. Fixing two points $x\not=y$, we have~:

$$\underset{N\rightarrow+\infty}{\lim}a_N(x,y)+\alpha(x,y)a_N(y,x)=\E_x\({\sum_{n=0}^{T_y-1}1_{S_n=x}}\)+1_{TR}G(x,x)P_x(T_y=\infty)P_y(T_x=\infty).$$
\end{lemme}

\noindent
{\it Proof of the lemma~:}

\noindent
Let again $c=\E_x(\sum_{0\leq n<T_y}1_{S_n=x})$ and $d=\E_y(\sum_{0\leq n<T_x}1_{S_n=y})$. Using \eqref{import} and \eqref{deb}~:

\begin{eqnarray}
da_N(x,y)+ca_N(y,x)&=&d(G_N(x,x)-G_N(y,x))+c(G_N(y,y)-G_N(x,y))\nonumber\\
&=&d\[{c_N+\sum_{k=1}^{N-1}P_x(T_y=k)(G_{N-k}(y,x)-G_N(y,x))-P_x(T_y\geq N)G_N(y,x)}\]\nonumber\\
&+&c\[{\sum_{k=1}^{N-1}P_x(T_y=k)(G_{N}(y,y)-G_{N-k}(y,y))+P_x(T_y\geq N)G_N(y,y)}\].\nonumber\end{eqnarray}

\medskip
\noindent
Set $b_N=cG_N(y,y)-dG_N(y,x)=cG_N(y,y)-dG_N(x,x)+d(G_N(x,x)-G_N(y,x))$. By the beginning of the section and the {\it claim} of the first section, $(b_N)$ is bounded. Therefore~:

$$da_N(x,y)+ca_N(y,x)=dc_N+\sum_{k=1}^{N-1}P_x(T_y=k)(b_N-b_{N-k})+P_x(T_y\geq N)b_N.$$

\medskip
\noindent
Let us study the limit of each term in the right-hand side, as $N\rightarrow+\infty$. The first one tends to $cd$. For the other ones, we distinguish the natural cases~:

\medskip
\noindent
- Transience. Then $b_N\rightarrow cG(y,y)-dG(y,x)=d(G(x,x)-G(y,x))$. By dominated convergence, the second term goes to zero. The limit thus exists and equals~:

$$cd+P_x(T_y=\infty)d(G(x,x)-G(y,x))=cd+P_x(T_y=\infty)P_y(T_x=\infty)dG(x,x).$$

\medskip
\noindent
- Null recurrence. Then $P_u(S_n=v)\rightarrow0$, for any $u$, $v$. This gives $b_N-b_{N+1}\rightarrow0$ and so $b_N-b_{N-k}\rightarrow 0$ for fixed $k$. By dominated convergence the second term goes to zero. As $P_x(T_y\geq N)b_N\rightarrow0$, the limit is thus $cd$ in this case.

\medskip
\noindent
- Positive recurrence. Again the third term tends to $0$. Aperiodicity implies that $P_u(S_n=v)\rightarrow\pi(v)$, where $\pi$ is the invariant probability measure for the chain. For fixed $k$, $b_N-b_{N-k}\rightarrow ck\pi(y)-dk\pi(x)=0$, as $\pi(x)/\pi(y)=c/d$ in this case. By dominated convergence once more the second term goes to $0$ and the limit also equals $cd$. 

\medskip
This concludes the proof of the lemma.

\fin

\medskip
\noindent
\begin{remark}
This lemma in the recurrent case is due to Chung, see Kemeny-Snell-Knapp \cite{KSK}, Theorem 9.7. Notice that the proof is somehow identical to that in {\it Step 2} of Proposition \ref{inega}. Again, the right-hand side of the formula is essentially $G(x,x)$, when $y$ goes to infinity. The idea now would be to understand the left-hand side with analytical tools. The quantity $\alpha(x,y)$ has to be analyzed closely. For a homogeneous random walk, $\alpha(x,y)=1$, since $P_0(T_0<T_x)=P_0(T_0<T_{-x})$; cf for example the end of the preliminary section. 
\end{remark}

\section{Fourier transform of probability measures on $\N^*$}

We study here $1/(1-\hat{\mu}_+)$ when $\mu_+$ is a probability measure on $\N^*$. The random walk with step $\mu_+$ is transient, hence $Re(1/(1-\hat{\mu}_+))\in L^1(0,2\pi)$. In fact $(1/(2\pi))\int_0^{2\pi}Re(1/(1-\hat{\mu}_+(t)))dt\leq1$, directly by relation \eqref{cf} and Fatou's lemma, as $G^{\mu_+}(0,0)=1$. The exact value would follow easily from the considerations of Proposition \ref{inega} combined with renewal theory. 

\medskip
We instead make an Herglotz type computation, using complex analysis. This allows to derive the renewal theorem directly from the Riemann-Lebesgue lemma. We next show that the Fourier coefficients of $Re(1/(1-\hat{\mu}_+))$ have an interesting probabilistic interpretation.

\begin{lemme}
\label{herglotz}

$ $

\noindent
Let $\mu_+$, with $\mbox{Supp}(\mu_+)\subset \N^*$ and $\mbox{gcd(Supp}(\mu_+))=1$. Then $t\longmapsto Re(1/(1-\hat{\mu}_+(t))$ is real, positive, even and in $L^1(0,2\pi)$, with~:

\begin{equation}
\label{cte}
\frac{1}{2\pi}\int_0^{2\pi}Re\({\frac{1}{1-\hat{\mu}_+(t)}}\)dt=1-\frac{1}{2E^{\mu_+}(X)}.
\end{equation}

\end{lemme}

\noindent
{\it Proof of the lemma~:}

\noindent
Introduce $f(z)=1/(1-\E^{\mu_+}(z^X))$, holomorphic in $\Delta=\{|z|<1\}$. For $0\leq r\leq 1$, the map $z\longmapsto Re(f(rz))$ is $>0$ and harmonic on $\Delta$. Thus, fixing $0<r<1$ and using the Poisson kernel~:

\begin{equation}
\label{ree}
Re(f(rz))=\frac{1}{2\pi}\int_0^{2\pi}Re\({\frac{e^{i\theta}+z}{e^{i\theta}-z}}\)Re(f(re^{i\theta}))d\theta,~z\in\Delta.
\end{equation}

\noindent
By holomorphic extension (using that $f(0)=1\in\R$)~:

\begin{equation}
\label{holol}
f(rz)=\frac{1}{2\pi}\int_0^{2\pi}\frac{e^{i\theta}+z}{e^{i\theta}-z}Re(f(re^{i\theta}))d\theta,~z\in\Delta.
\end{equation}

\medskip
\noindent
Taking $z=0$ in \eqref{ree}, we get $1=\frac{1}{2\pi}\int_0^{2\pi}Re(f(re^{i\theta}))d\theta$, so the positive measures $(\nu_r)_{0<r<1}$ on the torus $\R\backslash 2\pi\Z$ with density $\theta\longmapsto Re(f(re^{i\theta}))$ have constant mass $2\pi$. Let us take a cluster value $\nu$ of $\nu_r$, as $r\uparrow1$, for the weak-$*$ topology. We get from \eqref{holol}, fixing first $z\in\Delta$~:

$$f(z)=\frac{1}{2\pi}\int_0^{2\pi}\frac{e^{i\theta}+z}{e^{i\theta}-z}d\nu(\theta)=\frac{1}{2\pi}\int_0^{2\pi}\({1+2\sum_{n\geq1}z^ne^{-in\theta}}\)d\nu(\theta).$$

\noindent
Permuting the sum and the integral in the last expression, the Fourier coefficients of $\nu$ are uniquely determined by the development in series of $f$ around 0. Hence $\nu$ is unique and we conclude that $(\nu_r)$ converges to $\nu$, as $r\uparrow1$. We shall now determine this measure.

\medskip
First, when $\theta\in\R\backslash2\pi\Z$ is fixed, then $\hat{\mu}_+(\theta)\not=1$, so $\lim_{r\uparrow1}Re(f(re^{i\theta}))=Re(1/(1-\hat{\mu}_+(\theta)))$. Thus $\nu$ is locally $Re(1/(1-\hat{\mu}_+(\theta)))d\theta$, giving $\nu=Re(1/(1-\hat{\mu}_+(\theta)))d\theta+\alpha_{0}\delta_{0}$. Hence~:

$$\frac{1}{1-\E^{\mu_+}(z^X)}=\frac{1}{2\pi}\int_0^{2\pi}\({\frac{e^{i\theta}+z}{e^{i\theta}-z}}\)Re\({\frac{1}{1-\hat{\mu}_+(\theta)}}\)d\theta+\frac{\alpha_{0}}{2\pi}\({\frac{1+z}{1-z}}\),~z\in\Delta.$$

\noindent
To determine $\alpha_0$, take $z=e^{-u}$, with a real $u\downarrow0$, and multiply both sides by $1-e^{-u}$~:

$$\frac{1-e^{-u}}{1-\E^{\mu_+}(e^{-uX})}=\frac{1}{2\pi}\int_0^{2\pi}(1-e^{-u})\({\frac{e^{i\theta}+e^{-u}}{e^{i\theta}-e^{-u}}}\)Re\({\frac{1}{1-\hat{\mu}_+(\theta)}}\)d\theta+\frac{\alpha_{0}}{2\pi}(1+e^{-u}).$$

\noindent
The left-hand side goes to $1/E^{\mu_+}(X)$, monotonically. As $(1-e^{-u})\times(e^{i\theta}+e^{-u})/(e^{i\theta}-e^{-u})$ stays bounded by 2, the first term on the right-hand side tends to 0 by dominated convergence. Finally we get $\alpha_0/\pi=1/\E^{\mu_+}(X)$ and therefore the relation~:

$$\frac{1}{1-\E^{\mu_+}(z^X)}=\frac{1}{2\pi}\int_0^{2\pi}\({\frac{e^{i\theta}+z}{e^{i\theta}-z}}\)Re\({\frac{1}{1-\hat{\mu}_+(\theta)}}\)d\theta+\frac{1}{2\E^{\mu_+}(X)}\({\frac{1+z}{1-z}}\),~z\in\Delta.$$

\noindent
Expression \eqref{cte} is now given by $z=0$.

\fin

\begin{prop} 
\label{four}

$ $

\noindent
Let $\mu_+$, with $\mbox{Supp}(\mu_+)\subset \N^*$ and $\mbox{gcd(Supp}(\mu_+))=1$. 

\medskip
\noindent
i) For $x\geq1$~:

$$\frac{1}{\pi}\int_0^{2\pi}\cos(tx)Re\({\frac{1}{1-\hat{\mu}_+(t)}}\)dt=P_0^{\mu_+}(T_x<\infty)-\frac{1}{E^{\mu_+}(X)}.$$

\medskip
\noindent
As a result (renewal theorem, Erd\"os-Feller-Pollard \cite{efp}), $\lim_{x\rightarrow+\infty}P_0^{\mu_+}(T_x<\infty)=1/E^{\mu_+}(X)$.

\medskip
\noindent
ii) The function $t\longmapsto t/|1-\hat{\mu}_+(t)|$ belongs to $L^2(0,\pi)$. The function $t\longmapsto Im(1/(1-\hat{\mu}_+(t))$ is real and odd; it does not belong to $L^1(0,\pi)$, whenever $E^{\mu_+}(X)<\infty$. Also, for $x\geq1$~:

$$\frac{1}{\pi}\int_0^{2\pi}\sin(tx)Im\({\frac{1}{1-\hat{\mu}_+(t)}}\)dt=P_0^{\mu_+}(T_x<\infty).$$

\medskip
\noindent
iii) We have $|1-\hat{\mu}_+|^{-1}\in L^{\gamma}(0,2\pi)$, $0<\gamma<1$. Also $t^{\varepsilon}|1-\hat{\mu}_+(t)|^{-1}\in L^1(0,\pi)$, $\varepsilon>0$.

\end{prop}

\noindent
{\it Proof of the proposition~:}

\noindent
$i)$ Start as in Prop. \ref{inega}. Fixing $x\geq1$, we first have, for $N\geq1$~:

$$2G_N^{\mu_+}(0,0)-G_N^{\mu_+}(0,x)-G_N^{\mu_+}(0,-x)=\frac{1}{\pi}\int_0^{2\pi}\frac{1-\cos(tx)}{1-\hat{\mu}_+(t)}(1-(\hat{\mu}_+(t))^N)~dt.$$ 

\noindent
We can again take the limit as $N\rightarrow+\infty$ in the right-hand side and next the real part. The limit of the left-hand side is trivial, so for any $x\geq1$~:

\begin{equation}
\label{posit}
2-P_0^{\mu_+}(T_x<\infty)=\frac{1}{\pi}\int_0^{2\pi}(1-\cos (tx))Re\({\frac{1}{1-\hat{\mu}_+(t)}}\)~dt.
\end{equation}

\noindent
The fact that $Re(1/(1-\hat{\mu}_+))\in L^1(0,2\pi)$ can be recovered when minoring the right-hand side by $\pi^{-1}\int_{\delta}^{2\pi-\delta}$, $\delta>0$, letting $x\rightarrow+\infty$ with the Riemann-Lebesgue lemma and finally $\delta\rightarrow0$.

\medskip
\noindent
Subtracting \eqref{posit} to twice \eqref{cte}, we obtain the desired relation for $x\geq1$. Then the renewal theorem is a consequence of the Riemann-Lebesgue lemma. Mention in passing another proof, even simpler. By \eqref{posit} and the Riemann-Lebesgue lemma, $\lim_{x\rightarrow+\infty}P_0(T_x<\infty)$ exists; then in Spitzer \cite{spitz}, P3, the limit is identified as $1/\E^{\mu_+}(X)$.

\medskip
\noindent
$ii)$ Let us place on $(0,\pi)$. As $1-Re(\hat{\mu}_+(t))\geq\alpha t^2$, we have~:

$$\frac{\alpha t^2}{|1-\hat{\mu}_+(t)|^2}\leq \frac{1-Re(\hat{\mu}_+(t))}{|1-\hat{\mu}_+(t)|^2}=Re\({\frac{1}{1-\hat{\mu}_+(t)}}\)\in L^1(0,\pi).$$

\noindent
For the imaginary part, let us write, fixing $x\geq1$ and taking $N\geq1$~:

\begin{eqnarray}
G^{\mu_+}_N(0,x)-G^{\mu_+}_N(0,-x)&=&\frac{1}{2\pi}\int_0^{2\pi}(e^{-itx}-e^{itx})\sum_{n=0}^{N-1}(\hat{\mu}_+(t))^ndt\nonumber\\
&=&-\frac{i}{\pi}\int_0^{2\pi}\sin(tx)\({\frac{1-(\hat{\mu}_+(t))^N}{1-\hat{\mu}_+(t)}}\)dt.\end{eqnarray}

\noindent
As $\sin(tx)Im(1/(1-\hat{\mu}_+(t)))$ is integrable, we can let $N\rightarrow+\infty$ and then take the imaginary part inside the integral. The left-hand side limit being obvious, we obtain for $x\geq1$~:

$$P_0^{\mu_+}(T_x<\infty)=\frac{1}{\pi}\int_0^{2\pi}\sin(tx)Im\({\frac{1}{1-\hat{\mu}_+(t)}}\)dt.$$

\noindent
Whenever $E^{\mu_+}(X)<\infty$, the left-hand side goes to $1/E^{\mu_+}(X)>0$, hence the Riemann-Lebesgue lemma is not verified, giving $Im(1/(1-\hat{\mu}_+(t)))\not\in L^1(0,2\pi)$.

\medskip
\noindent
$iii)$ The holomorphic function $f(z)=1/(1-\E^{\mu_+}(z^X))$, $z\in \Delta$, has a positive and harmonic real part. The latter thus is in $h^1(\Delta)$. By Duren \cite{duren}, Theorem 4.2, $f\in H^{\gamma}(\Delta)$, $0<\gamma<1$, i.e.~:

$$\sup_{0<r<1}\int_0^{2\pi}|f(re^{i\theta})|^{\gamma}~d\theta<\infty.$$

\noindent
By Fatou's lemma, as $r\uparrow1$, we get $\int_0^{2\pi}|1-\hat{\mu}_+(\theta)|^{-\gamma}~d\theta<\infty$. For the last point, write~:

$$\frac{t^{\varepsilon}}{|1-\hat{\mu}_+(t)|}=\frac{t^{\varepsilon}}{|1-\hat{\mu}_+(t)|^{\varepsilon/2}}\times\frac{1}{|1-\hat{\mu}_+(t)|^{1-\varepsilon/2}}\in L^1(0,\pi),$$

\noindent
as the first term on the right-hand side is bounded. This ends the proof of the proposition. 

\fin

\medskip
\begin{cor}

\label{fourierr}

$ $

\noindent
$i)$ $Re((1-\hat{\mu}_+)^{-1})\in L^2(0,2\pi)$ iff $(G^{\mu_+}(0,x)-1/\E^{\mu_+}(X))_{x\geq0}\in l^2$.

\medskip
\noindent
$ii)$ $Im((1-\hat{\mu}_+)^{-1})\in L^2(0,2\pi)$ iff $(G^{\mu_+}(0,x))_{x\geq0}\in l^2$. In this case $|1-\hat{\mu}_+|^{-1}\in L^2(0,2\pi)$. 

\end{cor}

\noindent
{\it Proof of the corollary~:}

\noindent
Point $i)$ is clear as $Re((1-\hat{\mu}_+)^{-1})\in L^1(0,2\pi)$, so the $(G^{\mu_+}(0,x)-1/\E^{\mu_+}(X))$ are its Fourier coefficients. Idem, when $Im((1-\hat{\mu}_+)^{-1})\in L^2(0,2\pi)$ then the $(G^{\mu_+}(0,x))_{x\geq0}$ are the corresponding Fourier coefficients and thus belong to $l^2$. Reciprocally, if $(G^{\mu_+}(0,x))_{x\geq0}$ is $l^2$, define the $L^2$ odd function $f(t)=\sum_{x\geq1}G^{\mu_+}(0,x)\sin (xt)$. For all $x\in\Z$, we thus have~:

$$\int_0^{2\pi}\frac{\sin(tx)}{\sin t}\[{\sin t\({Im((1-\hat{\mu}_+(t))^{-1})-f(t)}\)}\]dt=0.$$

\noindent
The function inside the brackets belongs to $L^2$ and is even. Writing $\sin(t(1+x)+\sin(t(1-x))=2\sin(t)\cos(tx)$, for $x\geq0$, the latter is thus orthogonal to all $\cos(tx)$, $x\geq0$, hence equals zero a.-e.. Hence $Im((1-\hat{\mu}_+)^{-1})=f$, a.-e., and thus belongs to $L^2$. 

\medskip
Finally, when $Im((1-\hat{\mu}_+)^{-1})\in L^2(0,2\pi)\subset L^1(0,2\pi)$, then $\E^{\mu_+}(X)=\infty$ by Proposition \ref{four}. The conditions on Fourier coefficients in order to belong to $L^2(0,2\pi)$, for $Re((1-\hat{\mu}_+)^{-1})$ and $Im((1-\hat{\mu}_+)^{-1})$, are now identical.

\fin

\bigskip
\noindent
\begin{remark} Using $Re(1-\hat{\mu}_+(t))=2\E^{\mu_+}(\sin^2(tX/2))\geq(2/\pi^2)t^2\E^{\mu_+}(X^21_{X<\pi/t})$, we get~:

$$\int_0^{\varepsilon}\frac{t^2\E^{\mu_+}(X^21_{X<\pi/t})}{|1-\hat{\mu}_+(t)|^2}~dt<\infty.$$

\noindent
This is a little improvement of Prop. \ref{four} $ii)$ when $X\not\in L^2$. In Proposition \ref{four} the Fourier coefficients of $Re((1-\hat{\mu}_+)^{-1})$ and $Im((1-\hat{\mu}_+)^{-1})$ are probabilistic quantities. Those of $Re((1-\hat{\mu}_+)^{-1})$ exactly measure the error in the renewal theorem. Another question is whether $t^{1/2+\varepsilon}/|1-\hat{\mu}_+(t)|\in L^2(0,\pi)$, for $\varepsilon>0$. Also, $\sum_{x\geq0}G^{\mu_+}(0,x)=+\infty$, hence $(G^{\mu_+}(0,x))_{x\geq0}$ is never $l^1$. As detailed in the next section, it is $l^2$ iff some symmetric oscillating random walk on $\Z$ is transient.  \end{remark}

\bigskip
\noindent
\begin{remark}
For complex numbers $a$ and $b$, write $\langle a,b\rangle=Re(a\bar{b})$ for the real inner product of the vectors in $\R^2$ with affixes $a$ and $b$. As a corollary of Prop. \ref{four}, although $1/|1-\hat{\mu}_+|$ may not belong to $L^1(0,2\pi)$, we have $\langle (1-\hat{\mu}_+(t))^{-1},e^{itx}\rangle\in L^1(0,2\pi)$, for all $x\in\Z$, with~:

$$\left\{{\begin{array}{c}
\frac{1}{\pi}\int_0^{2\pi}\langle (1-\hat{\mu}_+(t))^{-1},e^{itx}\rangle~dt=2G^{\mu_+}(0,x)-1/E^{\mu_+}(X),~x\geq0,\\
\\
\frac{1}{\pi}\int_0^{2\pi}\langle (1-\hat{\mu}_+(t))^{-1},e^{-itx}\rangle~dt=-1/E^{\mu_+}(X),~x\geq1.\end{array}}\right.$$
\end{remark}

\noindent
Let $x\geq1$. Even if this is not always true, in general $Im(\hat{\mu}_+(t))\geq0$ for small $t>0$, so in this case $1/(1-\hat{\mu}_+(t))$ is in the first quadrant, as well as $e^{itx}$ and contrary to $e^{-itx}$. Hence it seems natural that the first integral above is larger than the second one.

\medskip
To conclude this section, we present a variation on Lemma \ref{herglotz}. 

\begin{lemme}
\label{varia}

$ $

\noindent
Let $\mu_+$, with $\mbox{Supp}(\mu_+)\subset \N^*$ and $\mbox{gcd(Supp}(\mu_+))=1$. Then~:

$$\frac{1}{2\pi}\int_0^{2\pi}Re\({\frac{1}{1-\hat{\mu}_+(t)}}\)\frac{1}{1+\[{\frac{\sin(t/2)}{\sinh(1/2)}}\]^2}dt=\frac{\tanh(1/2)}{1-\E^{\mu_+}(e^{-X})}-\frac{1}{2\E^{\mu_+}(X)}.$$
\end{lemme}

\noindent
{\it Proof of the lemma~:} 

\noindent
Introduce the homography $\rho(z)=(1-z)/(1+z)$, exchanging the open unit disk $\Delta$ and the half plane $Re>0$. Let $f(z)=1/(1-\E^{\mu_+}(e^{-\rho(z)X})),~z\in \Delta$. It is holomorphic in $\Delta$ and $Re(f)$ is $>0$ and harmonic on $\Delta$. By harmonicity at $z=0$, for $0<r<1$~:

$$1/(1-\E^{\mu_+}(e^{-X}))=\frac{1}{2\pi}\int_0^{2\pi}Re(f(re^{i\theta})d\theta.$$ 

\noindent
Proceeding as in Lemma \ref{herglotz}, we get~:

$$f(z)=\frac{1}{2\pi}\int_0^{2\pi}\frac{e^{i\theta}+z}{e^{i\theta}-z}d\nu(\theta),~z\in\Delta,$$

\noindent
where $\nu$ is the limit as $r\uparrow1$ of the positive measures $\nu_r$ on $\R\backslash 2\pi\Z$ with density $Re(f(re^{i\theta}))$. In order to detail $\nu$, note first that when $\theta\in(-\pi,\pi)$ is fixed, then~:

$$\lim_{r\uparrow1}\rho(re^{i\theta})=\frac{1-e^{i\theta}}{1+e^{i\theta}}=-i\tan(\theta/2).$$

\noindent
When also $\theta\not\in\{2\arctan(2k\pi),~k\in\Z\}$, then $\E^{\mu_+}(e^{i\tan(\theta/2)X})\not=1$ and $\nu$ is locally $g(\theta)d\theta$, with $g(\theta)=Re(1/(1-\E^{\mu_+}(e^{i\tan(\theta/2)X})))$. Hence $\nu$ decomposes as~:

$$\nu=Re(1/(1-\E^{\mu_+}(e^{i\tan(t/2)X})))dt+\alpha_{\pi}\delta_{\pi}+\sum_{k\in\Z}\alpha_k\delta_{\theta_k},$$

\noindent
where $\theta_k=2\arctan(2\pi k)$, for non-negative $\alpha_{\pi}$ and $(\alpha_k)$. In order to determine these coefficients, start from the relation, for $z\in\Delta$~:

$$\frac{1}{1-\E^{\mu_+}(e^{-\rho(z)X})}=\frac{1}{2\pi}\int_{-\pi}^{\pi}\({\frac{e^{i\theta}+z}{e^{i\theta}-z}}\)g(\theta)d\theta+\frac{\alpha_{\pi}}{2\pi}\({\frac{1-z}{1+z}}\)+\sum_{k\in\Z}\frac{\alpha_k}{2\pi}\({\frac{e^{i\theta_k}+z}{e^{i\theta_k}-z}}\).$$

\noindent
Recall that $g$ is integrable and $\sum_k\alpha_k<\infty$ (the mass of $\nu$ is $2\pi/(1-\E^{\mu_+}(e^{-X}))$). Take $z=-r$ above, as $r\uparrow1$, and multiply first both sides by $1-r$. Notice that for $\theta\in(-\pi,\pi)$, $(1-r)\times(e^{i\theta}-r)/(e^{i\theta}+r)$ stays bounded by 2 and converges to 0. Hence as $r\uparrow1$, by dominated convergence, the right-hand side converges to $\frac{\alpha_{\pi}}{2\pi}(1+1)$, whereas the left-hand side is equivalent to $(1-r)$ and therefore goes to 0. We obtain $\alpha_{\pi}=0$.

\medskip
\noindent
Fixing $k\in\Z$, take now $z=re^{i\theta_k}$ and let $r\uparrow1$, after multiplying both sides  by $(1-r)$. Idem, for $\theta\in(-\pi,\pi)\backslash\{\theta_k\}$, $(1-r)\times(e^{i\theta}+re^{i\theta_k})/(e^{i\theta}-re^{i\theta_k})$ stays bounded by 2 and converges to 0. By dominated convergence the right-hand side converges to $\alpha_k/\pi$, and this equals~:

$$\lim_{r\uparrow1}\frac{1-r}{1-\E^{\mu_+}(e^{-\rho(re^{i\theta_k})X})}.$$

To determine the limit, note first that $\lim_{u\downarrow 0^+}u/(1-\E^{\mu_+}(e^{-uX}))=1/\E^{\mu_+}(X)$, by monotone convergence. Next, $\rho(e^{i\theta_k})=-i\tan(\theta_k/2)=-2ik\pi$, so the denominator is~:

$$1-\E^{\mu_+}(e^{-(\rho(re^{i\theta_k})-\rho(e^{i\theta_k}))X}).$$

\noindent
We next have, decomposing in real and imaginary parts~:

\begin{eqnarray}
\rho(re^{i\theta_k})-\rho(e^{i\theta_k})=\frac{2(1-r)e^{i\theta_k}}{(1+re^{i\theta_k})(1+e^{i\theta_k})}&=&\frac{1-r}{\cos(\theta_k/2)}\frac{(1+r)\cos(\theta_k/2)+i(1-r)\sin(\theta_k/2)}{(1+r)^2\cos^2(\theta_k/2)+(1-r)^2\sin^2(\theta_k/2)}\nonumber\\
&=&A(r)+iB(r).\nonumber\end{eqnarray}

\noindent
- Case 1 : $\E^{\mu_+}(X)=+\infty$. Then, as $A(r)/(1-r)\rightarrow_{r\uparrow1}1/(2\cos^2(\theta_k/2))>0$~:

$$\left|{\frac{1-r}{1-\E^{\mu_+}(e^{-\rho(re^{i\theta_k})X})}}\right|\leq\frac{1-r}{1-\E^{\mu_+}(e^{-A(r)X})}\rightarrow_{r\uparrow1}0.$$

\noindent
- Case 2 : $\E^{\mu_+}(X)<+\infty$. Then~:

$$\frac{1-\E^{\mu_+}(e^{-(\rho(re^{i\theta_k})-\rho(e^{i\theta_k}))X})}{1-r}=\frac{1-\E^{\mu_+}(e^{-A(r)X})}{1-r}+\E^{\mu_+}(e^{-A(r)X}(1-e^{-iB(r)X})/(1-r)).$$

\noindent
As $r\uparrow1$, the first term on the right-hand side tends to $\E^{\mu_+}(X)/(2\cos^2(\theta_k/2))$. Since $t\longmapsto e^{it}$ is 1-Lipschitz on $\R$, $|1-e^{-iB(r)X})|/(1-r)\leq |B(r)|X/(1-r)$. As $|B(r)|\leq C(1-r)^2$, the previous quantity is both bounded by $C'X$ and tends to 0 as $r\uparrow1$. Since $\E^{\mu_+}(X)<\infty$, by dominated convergence the second term goes to 0 as $r\uparrow1$. Finally, $\alpha_k/\pi=2\cos^2(\theta_k/2)/\E^{\mu_+}(X)$. As $\cos^2(\theta_k/2)=1/(1+4\pi^2k^2)$, this leads to~:

\begin{eqnarray}
\frac{1}{1-\E^{\mu_+}(e^{-\rho(z)X})}&=&\frac{1}{2\pi}\int_{-\pi}^{\pi}\({\frac{e^{i\theta}+z}{e^{i\theta}-z}}\)Re\({\frac{1}{1-\E^{\mu_+}(e^{i\tan(\theta/2)X})}}\)d\theta\nonumber\\
&+&\frac{1}{\E^{\mu_+}(X)}\sum_{k\in\Z}\({\frac{e^{i\theta_k}+z}{e^{i\theta_k}-z}}\)\frac{1}{1+4\pi^2k^2},~z\in\Delta.\nonumber
\end{eqnarray}

\noindent
Taking $z=0$ and making the change of variable $\theta=2\arctan t$ in the first integral~:

\begin{eqnarray}
\frac{1}{1-\E^{\mu_+}(e^{-X})}&=&\frac{1}{\pi}\int_{-\infty}^{\infty}Re\({\frac{1}{1-\hat{\mu}_+(t)}}\)\frac{1}{1+t^2}dt+\frac{1}{\E^{\mu_+}(X)}\sum_{k\in\Z}\frac{1}{1+4\pi^2k^2}\nonumber\\
&=&\frac{1}{\pi}\int_{0}^{2\pi}Re\({\frac{1}{1-\hat{\mu}_+(t)}}\)\sum_{k\in\Z}\frac{1}{1+(t+2k\pi)^2}dt+\frac{1}{\E^{\mu_+}(X)}\sum_{k\in\Z}\frac{1}{1+4\pi^2k^2}.\nonumber\end{eqnarray}

\noindent
Finally, for a real $a$ and a complex number $z$, we have (cf Cartan \cite{cartan}, ex. 4, p172)~: 

$$\frac{\pi}{a}\frac{\sinh(2\pi a)}{\cosh(2\pi a)-\cos(2\pi z)}=\sum_{k\in\Z}\frac{1}{a^2+(z+k)^2}.$$

\noindent
Taking $z=t/(2\pi)$ and $a=1/(2\pi)$, we get $\sum_{k\in\Z}1/(1+(t+2k\pi)^2)=(1/2)\frac{\sinh(1)}{\cosh(1)-\cos t}$ and we arrive at the announced formula. 

\fin

\medskip
\noindent
\begin{remark} In the definition of $z\longmapsto 1/(1-\E^{\mu_+}(e^{-\rho(z)X}))$, $z\in\Delta$, the term $\E^{\mu_+}(e^{-\rho(z)X})$ could be replaced by $\E^{\mu_+}(h(z)^X)$ or $\E^{\mu_+}(h(z^X))$, for any $h$ holomorphic in $\Delta$ with $|h|<1$ in $\Delta$, for example an automorphism of $\Delta$. One would get new relations, but concretely $\lim_{r\uparrow1}\E^{\mu_+}(h(re^{i\theta})^X)$ is often delicate to determine.

\end{remark}

\section{Kemperman's oscillating random walk}

Taking $\Z$ as state space and considering an inhomogeneous model, we now focus on oscillating random walks, as introduced by Kemperman in \cite{kemperman}. Fix a Markov chain, written as $(S_n)$, which jumps according to probability measures $\mu$ on $(-\infty,0]$ and $\nu$ on $[1,+\infty)$, respectively. In view of later applications, no moment assumption is made on either $\mu$ or $\nu$.

\medskip
Associate to each $\mu$ and $\nu$ one-sided sub-probability measures $\mu_+$ and $\nu_-$, respectively, with $\mbox{Supp}(\mu_+)\subset\N^*$ and $\mbox{Supp}(\nu_-)\subset-\N^*$, in the following way. For any $k\geq1$~:

\begin{equation}
\label{defmup}
\mu_+(k)=\P_0^{\mu}(S_n \mbox{ enters }\N^*\mbox{ at }k)\mbox{ and }\nu_-(-k)=\P_0^{\nu}(S_n\mbox{ enters }-\N^*\mbox{ at }-k).
\end{equation}

\smallskip
\noindent
The question of the recurrence of $0$ for $(S_n)$ will appear as a property of the sole couple $(\mu_+,\nu_-)$. We first discuss the relation between $\mu$ and $\mu_+$ (or $\nu$ and $\nu_-$). Classically, $\mu_+$ is related to the right Wiener-Hopf factor of $\mu$. The recurrence problem involves some inner product of the right Wiener-Hopf factor of $\mu$ with the left Wiener-Hopf factor of $\nu$.

\subsection{Link between $\mu$ and $\mu_+$}

For a measure $w$ on $\Z$, define the restrictions $w^-=w1_{\leq 0}$ and $w^+=w1_{\geq1}$. We place in the commutative Banach algebra of signed measures on $\Z$, with convolution as product, written as $w_1w_2$. Recall the fundamental property of the exponential, $\exp(w_1+w_2)=\exp(w_1)\exp(w_2)$, as well as the following identity for a non-negative measure $w$ with mass $<1$~:

$$\delta_0-w=\exp(-L_w)\mbox{, where }L_w=\sum_{n\geq1}\frac{w^n}{n}.$$

\noindent
Given a probability measure $w$ on $\Z$, write $(S_n^{w})$ for the $i.i.d.$ random walk with step $w$, with $S_0^w=0$. When several $(S_n^{w})$ appear, corresponding to different probability measures, they are supposed to be independent.

\begin{prop}

$ $

\noindent
Let $\mu$ be a probability measure on $\Z$ and $\mu_+$ defined as in \eqref{defmup}. Then $\mu_+$ is a probability measure iff $\sum_{n\geq1}\mu^n(\N^*)/n=+\infty$. When $\hat{\mu}_+(t)\not=1$~:

$$\frac{1}{1-\hat{\mu}_+(t)}=\lim_{s\uparrow1}e^{\sum_{n\geq1}s^n\widehat{(\mu^n)^+}(t)/n}.$$

\end{prop}

\noindent
{\it Proof of the proposition~:}

\noindent
Let $0<s<1$ and define $L_{\mu}^{\pm}=\sum_{n\geq1}s^n(\mu^n)^{\pm}/n$. Then $\delta_0-s\mu=\exp(-L_{\mu}^+)\exp(-L_{\mu}^-)$. Set $N=\min\{n\geq1,~S^{\mu}_n\geq1\}$, $\eta_0=\delta_0$ and $\eta_n(A)=P^{\mu}(N\geq n,~S_n\in A)$, $n\geq 1$. Let $\eta=\sum_{n\geq0}s^n\eta_n$. 

\medskip
\noindent
By definition, $\eta_{n+1}=(\eta_n)^-\mu$. Summing on $n\geq0$ with coefficients $s^{n+1}$, we get $\eta-\delta_0=\eta^-s\mu$. This gives $\eta^-(\delta_0-s\mu)=\delta_0-\eta^+$ and therefore $\eta^-\exp(-L_{\mu}^-)=(\delta_0-\eta^+)\exp(L_{\mu}^+)$. The left-hand side is a measure on $(-\infty,0]$ and the right-hand side on $[0,+\infty)$, with mass at 0 equal to one.

\medskip
\noindent
Hence $\eta^-\exp(-L_{\mu}^-)=(\delta_0-\eta^+)\exp(L_{\mu}^+)=\delta_0$. This gives $\delta_0-\eta^+=\exp(-L_{\mu}^+)$ or equivalently $\sum_{n\geq0}(\eta^+)^n=\exp(L_{\mu}^+)$, from which the assertions follow (observing that $\mu_+=\lim_{s\uparrow1}\eta^+$). 

\fin

\bigskip

\noindent
\begin{remark} Can we write the limit as $e^{\sum_{n\geq1}\widehat{(\mu^n)^+}(t)/n}$ for $0<t<2\pi$, when $\mu_+$ is a probability and $|\widehat{\mu}|<1$ on $(0,2\pi)$, hence suppressing the disgracious $\lim_{s\uparrow 1}$ ?
\end{remark}

\subsection{The concentrated Markov chain}

\begin{lemme}

$ $

\noindent
i) If either $\mu_+(\N^*)<1$ or $\nu_-(-\N^*)<1$, then $(S_n)$ is transient. 

\medskip
\noindent
ii) Let $\mu_+(\N^*)=\nu_-(-\N^*)=1$ and call $(Z_n)$ the Markov chain jumping with $\mu_+$ on $(-\infty,0]$ and $\nu_-$ on $[1,+\infty)$. Then $0$ is recurrent for $(S_n)$ is recurrent if and only if $0$ is recurrent for $(Z_n)$.
\end{lemme}

\noindent
{\it Proof of the lemma~:}

\noindent
$i)$ If $\mu_+(\N)<1$, then $(S_n^{\mu})$ a.-s. makes only finitely many records in the right direction, giving $S_n^{\mu}\rightarrow-\infty$, a.-s.. Thus $P_0(S_n^{\mu}\rightarrow-\infty\mbox{, with }S^{\mu}_k\leq 0,~\forall k\geq0)>0$ and so $P_0(S_n\rightarrow-\infty\mbox{ and }S_k\leq 0,~\forall k\geq0)>0$. Hence $(S_n)$ is transient. The situation $\nu_-(-\N^*)<1$ is treated similarly.

\medskip
\noindent
$ii)$ Let $\mu_+(\N^*)=\nu_-(-\N^*)=1$. Then $(S_n)$ visits both $(-\infty,0]$ and $[1,+\infty)$ infinitely often, a.-s.. Start $(Z_n)$ at 0. Idem, start $(S_n)$ at 0 and let $\tau $ be its a.-s. finite entrance time in $[1,+\infty)$. Then $S_{\tau}$ has the law of $Z_1$. Looking now at $(S_{\tau+n})_{n\geq0}$ at left record times on $[1,+\infty)$ and right record times on $(-\infty,0]$, then $(S_{\tau+n})_{n\geq0}$ a.-s. comes back to 0 iff $(Z_n)_{n\geq1}$ does. \fin

\medskip
We now suppose that $\mu_+(\N^*)=\nu_-(-\N^*)=1$. This property may not be sufficient for recurrence (rarely, very large jumps across 0 may occur, ensuring $|S_n|\rightarrow+\infty$). Using the previous lemma, we focus on $(Z_n)$. The latter random walk is rather particular as it can essentially be reduced to the two sequences of positive and negative jumps (written in the order they appear).

\begin{lemme}
\label{calcul}
$ $

\noindent
i) Let $0\leq k\leq n$ and $x,y\in\Z$. The sequences $(l_i^+)_{1\leq i\leq k}$ and $(-l_j^-)_{1\leq j\leq n-k}$ are respectively the ordered sequences of positive and negative jumps of a trajectory (which is then unique) of $(Z_m)_{0\leq m\leq n}$, with $Z_0=x$ and $Z_n=y$, if and only if $\sum_{1\leq i\leq k}l_i^+-\sum_{1\leq j\leq n-k}l_j^-=y-x$ and $(l^+_k\geq y\mbox{, if }k\geq1~\&~ y>1)$ and $(-l^-_{n-k}\leq y-1\mbox{, if }n-k\geq 1~\&~ y\leq -1)$.

\medskip
\noindent
ii) We have $P_x(Z_n=y)=\sum_{k=0}^nP\({S_k^{\mu_+}+S_{n-k}^{\nu_-}=y-x,~X_k^{\mu_+}\geq y,~X^{\nu_-}_{n-k}\leq y-1}\)$, for any $n\geq0$ and $x,y\in\Z$. The second and third conditions disappear if $k=0$ or $n-k=0$. 

\end{lemme}

\noindent
{\it Proof of the lemma~:}

\noindent
$i)$ Starting from a trajectory of $(Z_m)_{0\leq m\leq n}$, denote by $(l_i^+)_{1\leq i\leq k}$ and $(-l_j^-)_{1\leq j\leq n-k}$, for some $0\leq k\leq n$, the ordered sequences of positive and negative jumps. From these two sequences we recover the trajectory, simply noting that the current position of the walker gives the direction of the next jump. Hence, starting for example from $x\leq 0$, use first the $(l_i^+)$ until reaching $[1,+\infty)$, next the $(l_j^-)$ until coming back to $(-\infty,0]$, etc, until exhausting the two lists. 

\medskip
Starting from $x$ and arriving at $y$, we have $\sum_i l_i^+-\sum_j l_j^-=y-x$. Suppose that $k\geq1$ and $y>1$ (the cases $n-k\geq1$ and $y\leq -1$ would be treated in the same way). When running the exhaustion process of the lists, two cases may occur~:

\noindent
- the $(l^+_i)$ are finished first. When this happens, the position is $>y$ and the last positive jump must have been $>y$. The path to $y$ is ended with the remaining negative jumps. 

\noindent
- the list $(l_i^-)$ is ended first. The trajectory terminates with positive jumps (each with a starting point in $(-\infty,0]$) and the last one has to be $\geq y$.

\medskip
Reciprocally, suppose the conditions satisfied and for instance $k\geq1~\&~y\geq 1$. Starting from $x$, run the exhaustion process of the lists. If the $(l^-_j)$ finish first, only positive jumps remain. As the last one is $\geq y$, this last sequence of jumps will have non-positive starting points, so the trajectory will be ``admissible". If the $(l^+_i)$ are ended first, we are $>y$ when this happens. Only remain negative jumps for going to $y$, hence the trajectory is also ``admissible".

\medskip
\noindent
$ii)$ Let independent $(X_k^{\mu_+},X_l^{\nu_-})_{k,l\geq0}$, with ${\cal L}(X_k^{\mu_+})=\mu_+$, ${\cal L}(X_l^{\nu_-})=\nu_-$. When running from some fixed $x$ the exhaustion process with the two lists $(X_k^{\mu_+})_{k\geq0}$ and $(X_l^{\nu_-})_{l\geq0}$, we obtain a realization of $(Z_n)_{n\geq0}$, with $Z_0=x$. By $i)$, $\{Z_n=y\mbox{, with }k \mbox{ positive jumps}\}=\{S_k^{\mu_+}+S_{n-k}^{\nu_-}=y-x,~X_k^{\mu_+}\geq y,~X^{\nu_-}_{n-k}\leq y-1\}$. Take the probability and sum on $0\leq k\leq n$ to get the result.

\fin

\medskip
The Green function of $(Z_n)$ is related to some $l^2$-inner product of the Green functions $G^{\mu_+}$ and $G^{\nu_-}$. See Theorems 4.6 and 4.8 of Kemperman \cite{kemperman}.

\begin{prop} 
\label{ggreen}

$ $

\noindent
i) The Green function $G^Z$ of $(Z_n)$ verifies $G^Z(0,0)=\sum_{m\geq0}G^{\mu_+}(0,m)G^{\nu_-}(0,-m)$ and~:

$$G^Z(0,0)=\lim_{s\uparrow 1}\frac{1}{2\pi}\int_0^{2\pi}Re\({(1-s\hat{\mu}_+(t))^{-1}(1-s\hat{\nu}_-(t))^{-1} }\)dt.$$

\noindent
ii) When $\nu_-(A)=\mu_+(-A)$, for $A\subset\Z$, then $G^Z(0,0)=(2\pi)^{-1}\int_0^{2\pi}|1-\hat{\mu}_+(t)|^{-2}dt$.

\end{prop}

\medskip
\noindent
\textit{Proof of the proposition~:}

\noindent

\noindent
$i)$ By Lemma \ref{calcul}, $P_0(Z_n=0)=\sum_{k=0}^nP\({S_k^{\mu_+}+S_{n-k}^{\nu_-}=0}\)=\sum_{k=0}^n(\mu_+^k\nu_-^{n-k})(0)$. Thus~:

$$G^Z(0,0)=\sum_{k,l\geq0}(\mu_+^k\nu_-^{l})(0)=\sum_{m\geq0}\sum_{k,l\geq0} \mu_+^k(m)\nu_-^l(-m)=\sum_{m\geq0}G^{\mu_+}(0,m)G^{\nu_-}(0,-m).$$

\noindent
Taking now $0<s<1$~:

$$\sum_{k,l\geq0}s^{k+l}(\mu_+^k\nu_-^{l})(0)=\frac{1}{2\pi}\int_0^{2\pi}\sum_{k,l\geq0}(s^k\hat{\mu}_+^ks^l\hat{\nu}_-^l)(t)dt=
\frac{1}{2\pi}\int_0^{2\pi}Re\({\frac{1}{1-s\hat{\mu}_+(t)}\frac{1}{1-s\hat{\nu}_-(t)}}\)dt.$$

\noindent
As $G^Z(0,0)=\lim_{s\uparrow 1}\sum_{k,l\geq0}s^{k+l}(\mu_+^k\nu_-^{l})(0)$, we get the desired expression.

\medskip
\noindent
$ii)$ In this case, $\hat{\mu}_+(t)$ is the conjugate of $\hat{\nu}_-(t)$, $t\in\R$. By $i)$, $2\pi G^Z(0,0)$ equals~:

$$\lim_{s\uparrow 1}\int_0^{2\pi}|1-s\hat{\mu}_+(t)|^{-2}dt=\lim_{s\uparrow 1}s^{-2}\int_0^{2\pi}|s^{-1}-\hat{\mu}_+(t)|^{-2}dt=\int_0^{2\pi}|1-\hat{\mu}_+(t)|^{-2}dt,$$

\noindent
where monotone convergence is used at the end (this does not work in general). 

\fin

\medskip

\noindent
\begin{remark}
The proof of $i)$ also gives the following interesting relation, interpreting of the recurrence criterion in terms of intersections of two independent random walks~:

$$G^Z(0,0)=\E\({\sum_{k,l\geq0}1_{S_k^{\mu_+}=-S_l^{\nu_-}}}\)=\E(\mbox{card}(\{S_k^{\mu_+},k\geq0\}\cap\{-S_l^{\nu_-},l\geq0\})).$$
\end{remark}


In fact $(Z_n)$ admits a natural invariant measure and this provides a characterization of positive recurrence. This was independently shown by Vo in \cite{vo}, which also exhibits an invariant measures for a related process, leading to an interesting recurrence condition.

\begin{prop}
\label{recpo}

$ $

\noindent
Let $\mu_+(\N^*)=\nu_-(-\N^*)=1$. Suppose that $(Z_n)$ is irreducible.

\medskip
\noindent
i) The measure $\pi(y)=\mu_+(\geq y)1_{y\geq1}+\nu_-(\leq y-1)1_{y\leq0}$ is invariant for $(Z_n)$. Hence $(Z_n)$ is positive recurrent iff $\E^{\mu_+}(X)<\infty$ and $\E^{\nu_-}(X)>-\infty$.

\medskip
\noindent
ii) If among $\E^{\mu_+}(X)$ and $\E^{\nu_-}(X)$ exactly one is finite, then $(Z_n)$ is null recurrent. 

\medskip
\noindent
iii) If $\int_0^{2\pi}|1-\hat{\mu}_+(t)|^{-2}dt<\infty$ and $\int_0^{2\pi}|1-\hat{\nu}_-(t)|^{-2}dt<\infty$, then $(Z_n)$ is transient, hence $(S_n)$.
\end{prop}

\noindent
{\it Proof of the proposition~:}

\noindent
Let the measure $\pi(y)=\mu_+(\geq y)1_{y\geq1}+\nu_-(\leq y-1)1_{y\leq0}$. Taking $y_0\geq1$, we have~:

\begin{eqnarray}
&&\sum_{y>y_0}\mu_+(\geq y)\nu_-(y_0-y)+\sum_{y\leq 0}\nu_-(\leq y-1)\mu_+(y_0-y)\nonumber\\
&=&\sum_{z>y_0}\mu_+(z)\sum_{y_0<y\leq z}\nu_-(y_0-y)+\sum_{z\geq y_0}\mu_+(z)\nu_-(\leq y_0-z-1)\nonumber\\
&=&\sum_{z>y_0}\mu_+(z)(\nu_-([y_0-z,-1])+\nu_-(\leq y_0-z-1))+\mu_+(y_0)\nu_-(\leq-1)=\mu_+(\geq y_0).\nonumber\end{eqnarray}

\noindent
The case $y_0\leq 0$ would be treated similarly. For $ii)$, if $\E^{\mu_+}(X)<\infty$, then $\lim_{m\rightarrow+\infty}G^{\mu_+}(0,m)=1/\E^{\mu_+}(X)>0$, by the renewal theorem. By Prop. \ref{ggreen}, $(Z_n)$ is recurrent iff $\sum_{m\geq0}G^{\nu_-}(0,-m)=+\infty$, which is true. The case $\E^{\nu_-}(X)>-\infty$ is similar.

\medskip
Concerning $ii)$, Prop. \ref{ggreen} gives that $(G^{\mu_+}(0,m))_{m\geq0}$ and $(G^{\nu_-}(0,-m))_{m\geq0}$ are in $l^2$, hence $(G^{\mu_+}(0,m)G^{\nu_-}(0,-m))_{m\geq0}$ is $l^1$, by the Cauchy-Schwarz inequality. Prop. \ref{ggreen} $i)$ then gives transience. This concludes the proof of the proposition.\fin

\subsection{Discussion}

1) A natural question now is whether $\int_0^{2\pi}Re((1-\hat{\mu}_+(t))^{-1}(1-\hat{\nu}_-(t))^{-1})dt$ always has a meaning, in the sense that $(Re((1-\hat{\mu}_+)^{-1}(1-\hat{\nu}_-)^{-1}))^-\in L^1(0,2\pi)$. Conditioned by a positive answer, the problem next is whether $G^Z(0,0)$ has the same order as $\int_0^{2\pi}Re((1-\hat{\mu}_+(t))^{-1}(1-\hat{\nu}_-(t))^{-1})dt$. A natural way seems to use Lemma \ref{fonda}.

\medskip
2) Starting from Prop. \ref{recpo}, the remaining case is $\E^{\mu_+}(X)=\E^{\nu_-}(-X)=+\infty$, with either $|1-\hat{\mu}_+(t)|^{-1}$ or $|1-\hat{\nu}_-(t)|^{-1}$ not in $L^2(0,2\pi)$. By Prop. \ref{four}, $Re((1-\hat{\mu}_+)^{-1})$ and $Re((1-\hat{\nu}_-)^{-1})$ have non-negative real cosine coefficients. Taking real functions $f>0$ and $g>0$ in $L^1(0,2\pi)$, even and with real non-negative Fourier coefficients $\hat{f}(n)$, $\hat{g}(n)$, $n\geq0$, a general question is whether $\int_0^{2\pi}f(t)g(t)dt\mbox{ and }\sum_{n\geq0}\hat{f}(n)\hat{g}(n)$ are infinite at the same time. We show the direction $``\leq"$, making a standard computations with the Fejer kernel $K_n$. We have~:

$$(K_n* f)(t)=\sum_{j=-n}^n\({1-\frac{|j|}{n+1}}\)\hat{f}(j)e^{ijt}.$$

\noindent
For $M>0$, $(1/2\pi)\int_0^{2\pi}(K_n* f)(g\wedge M)dt\leq \frac{1}{2\pi}\int_0^{2\pi}(K_n* f)gdt=\sum_{j=-n}^n(1-\frac{|j|}{n+1})\hat{f}(j)\hat{g}(j)$. For the left-hand side use that $K_n*f\rightarrow f$ in $L^1$. Using monotone convergence for the right-hand side, we obtain $(1/2\pi)\int_0^{2\pi}f(g\wedge M)dt\leq\sum_{n\geq0}\hat{f}(n)\hat{g}(n)$. Letting $M\rightarrow+\infty$, we get the result.

\medskip
\providecommand{\bysame}{\leavevmode\hbox to3em{\hrulefill}\thinspace}

\bigskip
{\small{\sc{Univ Paris Est Creteil, CNRS, LAMA, F-94010 Creteil, France\\
Univ Gustave Eiffel, LAMA, F-77447 Marne-la-Vall\'ee, France}}}

\it{E-mail address~:} {\sf julien.bremont@u-pec.fr}

\end{document}